\documentclass[12pt]{article}

 \usepackage{a4}                          
 \usepackage{amsmath}
 \usepackage{amssymb}
 \usepackage{graphicx}
 \newtheorem{thm}{Theorem}[section]

 \newtheorem{cor}[thm]{Corollary}
 \newtheorem{prop}[thm]{Proposition}

 \newtheorem{lemma}[thm]{Lemma}
 \newtheorem{rem}[thm]{Remark}
 \DeclareMathOperator{\Log}{Log}
 \DeclareMathOperator{\dilog}{dilog}
 \DeclareMathOperator{\Arg}{Arg}
\renewcommand{\t}{\theta}
 \newcommand{\eps}{\varepsilon}
 
\newcommand{\Int}{\int_{-\infty}^\infty}

 \title{A completely monotonic function used in an inequality of Alzer\footnote{Both authors acknowledge support by grant
10-083122 from The Danish Council for Independent Research $|$ Natural Sciences.}}

 \author{Christian Berg $^{\dagger}$ and Henrik L. Pedersen
   $^{\ddagger}$\\
 \footnotesize $\dagger$ \ Institute of Mathematical Sciences,
   University of Copenhagen\\
    \footnotesize Universitetsparken 5; DK-2100 K\o benhavn \O, Denmark\\
\footnotesize E-mail berg@math.ku.dk 
\\  \footnotesize
   $\ddagger$ \ Department of Basic Sciences and Environment\\
\footnotesize Faculty of Life Sciences, University of Copenhagen\\
\footnotesize Thorvaldsensvej 40, DK-1871 Frederiksberg C\\
\footnotesize E-mail henrikp@life.ku.dk
}

 \date{\today}

\begin{document}

 \maketitle

 \begin{abstract}
 The function $G(x)=\left(1-\ln x /\ln(1+x)\right)x\ln x$ has been
 considered by Alzer, Qi and Guo. We prove that $G'$ is completely
 monotonic by finding an integral representation of the holomorphic
 extension of $G$ to the cut plane. A main difficulty is caused by the
 fact that $G'$ is not a Stieltjes function. 
 \end{abstract}



 \section{Introduction and results}
In a recent paper \cite{Al}, Alzer proved a number of inequalities involving the
volume  of the unit ball in $\mathbb R^n$,
\begin{equation}\label{eq:vol}
\Omega_n=\frac{\pi^{n/2}}{\Gamma(1+n/2)},\quad n=1,2,\ldots.
\end{equation}
That paper contains many references to earlier results about
$\Omega_n$. We mention in particular that
Anderson and Qiu \cite{A:Q} proved that the sequence
$f(n)=\Omega_n^{1/(n\log n)},\;n\ge 2$ is strictly decreasing and
converges to $e^{-1/2}$. It is therefore of interest to study the
function
\begin{equation}\label{eq:f}
f(x)=\left(\frac{\pi^{x/2}}{\Gamma(1+x/2)}\right)^{1/(x\ln x)},
\end{equation}
and in \cite{B:P3} the authors have given an integral representation of $\log
f(x+1), x>0$ by considering its holomorphic extension to the cut plane
$\mathcal A=\mathbb C\setminus(-\infty,0]$.
From this representation it has been possible to deduce that $f(n+2)$ is a
Hausdorff moment sequence, in particular decreasing and convex.

The papers \cite{A:Q} and \cite{A:V:V} have also been an inspiration for several papers
about the functions
\begin{equation}\label{eq:Fa}
 F_a(x)=\frac{\ln \Gamma (x+1)}{x\ln(ax)}, \quad x>0, a>0,
 \end{equation}
see
\cite{Al},\cite{B:P1},\cite{B:P2},\cite{B:P3},\cite{E:L},\cite{Q:G}.
In particular, \cite{B:P3} contains an integral representation of the
meromorphic extension of $F_a$ to $\mathcal A$. From this representation it is
possible to deduce that $F_a$ is a Pick function if and only if $a\ge
1$. The relation between $F_a$ and $f$ is given by
$$
\log f(z+1)=\frac{\ln\sqrt{\pi}}{\Log(z+1)}-\frac12 F_2\left(\frac{z+1}2\right).
$$

 Alzer found the best constants $a^*,b^*$
such that for all $n\ge 2$
\begin{equation}\label{eq:alz}
\exp\left(\frac{a^*}{n(\log n)^2}\right)\le f(n)/f(n+1)<
\exp\left(\frac{b^*}{n(\log n)^2}\right).
\end{equation}

In the proof of this result Alzer considered
the function
\begin{equation}\label{eq:Greal}
G(x)=\left(1-\frac{\ln x}{\ln(1+x)}\right)x\ln x,
\end{equation}
and in \cite[Lemma 2.3]{Al} it was proved that $2/3<G(x)<1$ for
$x\ge 3$. Qi and Guo observed in
\cite{Q:G} that
$G$ is strictly increasing on $(0,\infty)$ with
$G((0,\infty))=(-\infty,1)$ and that $G(3)>2/3$, which gave a new
proof of the inequality  $2/3<G(x)<1$. Furthermore,  
in \cite[Remark 4]{Q:G} it was conjectured that 
\begin{equation}\label{eq:conj}
(-1)^{k-1}G^{(k)}(x)>0\;\; \mbox{for}\;\; x>0,k=1,2,\ldots,
\end{equation}
or equivalently that $G'$ is a completely monotonic function.

The main goal of this paper is to prove this conjecture. We do this by considering $G$ as a holomorphic function in the cut
plane $\mathcal A$. We put 
\begin{equation}\label{eq:G}
G(z)=\left(1-\frac{\Log z}{\Log(1+z)}\right)z\Log z,
\end{equation}
where $\Log z=\ln|z|+i \Arg z$ is the
principal logarithm in $\mathcal A$ and $-\pi<\Arg z<\pi$ for
$z\in\mathcal A$.  

Using the same Cauchy integral formula technique as in the paper
\cite{B:P3}, we shall establish the following theorem.

\begin{figure}[h]
\begin{center}
\includegraphics[width=100mm]{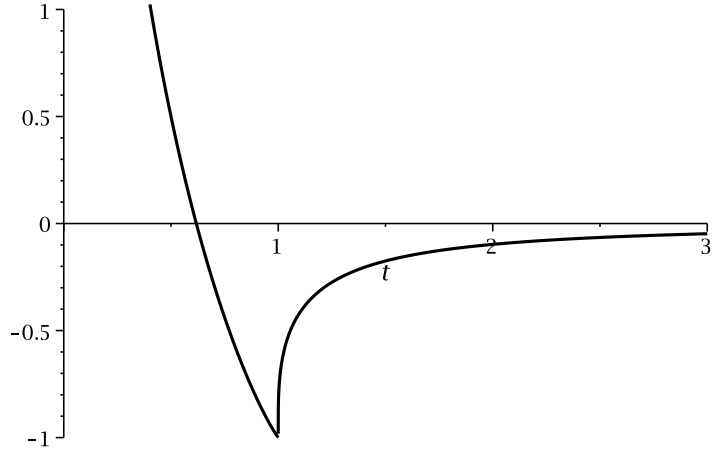}
\caption{The graph of $\rho$}
\label{fig:rho}
\end{center}
\end{figure}

\begin{thm}\label{thm:G1} The function $G$ from \eqref{eq:G} has the representation
\begin{equation}\label{eq:Grep}
 G(z)=1-\int_0^\infty\frac{\rho(t)}{z+t}\,dt,\quad z\in\mathcal A,
\end{equation}
where
\begin{equation}\label{eq:rho}
\rho(t)=\left\{\begin{array}{ll}
-\frac{t\ln ((1-t)/t^2)}{\ln(1-t)}\;, & \;\;\mbox{if}\;\; 0<t<1,\\

-\frac{t(\ln ((t-1)/t))^2}{(\ln(t-1))^2 + \pi^2}\;, & \;\;\mbox{if}\;\; 1<t<\infty.
\end{array}\right.
\end{equation}
\end{thm}
Notice that $\rho(1^-)=\rho(1^+)=-1$ so that $\rho$ is continuous on
the positive half-line. It is decreasing from  $\infty$ to $-1$ on the
interval $(0,1)$ with $\rho'(1^-)=-1$, and
increasing from $-1$ to $0$ on the interval $(1,\infty)$ with
$\rho'(1^+)=\infty$. We have 
$\rho((\sqrt{5}-1)/2)=0$. Notice also  that $\rho$ is integrable over
$(0,\infty)$ because of the asymptotics
$$
\rho(t)\sim -2\ln t \;\;\mbox{for}\;\; t\to 0^+;\quad \rho(t)\sim -\frac{1}{t(\ln
  t)^2}  \;\;\mbox{for}\;\;  t\to\infty.
$$
The graph of $\rho$ is shown in Figure \ref{fig:rho}.

Since $\rho$ assumes positive and negative
values, $G$ as well as $1-G$ are not  Stieltjes functions, but
nevertheless $1-G$ turns our to be completely monotonic, because
it is the Laplace transform of a positive function, as described in
the following theorem. In particular, $G'$ is completely monotonic so
 \eqref{eq:conj} holds. For
properties about completely monotonic functions and Stieltjes
functions we refer to \cite{B:F} and \cite{W}.

\begin{thm}\label{thm:G2} For $\Re z>0$ the function $1-G$ has the
  representation 
\begin{equation}\label{eq:1-G}
1-G(z)=\int_0^\infty e^{-zs}\varphi(s)\,ds,
\end{equation}
where
\begin{equation}\label{eq:phi}
\varphi(s)=\int_0^\infty e^{-st}\rho(t)\,dt>0 \;\;\mbox{for}\;\; s \ge
0.
\end{equation}
\end{thm}

\begin{figure}[h]
\begin{center}
\includegraphics[width=100mm]{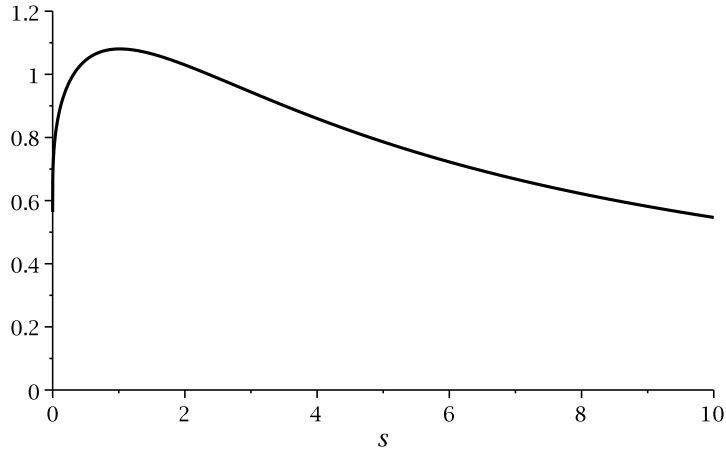}
\caption{The graph of $\varphi$}
\label{fig:phi}
\end{center}
\end{figure}

The graph of $\varphi$ is given in Figure \ref{fig:phi}.

The function $\varphi$ given in \eqref{eq:phi} is continuous and
bounded on $[0,\infty)$, but it is not integrable because
$1-G(x)\to\infty$ for $x\to 0^+$.

Setting $z=a+it$ in \eqref{eq:1-G} with $a>0$ we get:

\begin{cor}\label{thm:Bochner} (i) For each $a>0$
\begin{equation}\label{eq:L2}
1-G(a+it)=\int_0^\infty e^{-its}e^{-as}\varphi(s)\,ds,\quad
t\in\mathbb R
\end{equation}
is an analytic positive definite function of $t$, and it is the Fourier
transform of 
\begin{equation}\label{eq:plan}
e^{-as}\varphi(s)1_{[0,\infty)}(s).
\end{equation}

(ii) $G(a+it)-G(a)$ is a continuous negative definite function of $t$
for each $a>0$. In particular
\begin{equation}\label{eq:negdef}
\Re\,G(a+it)\ge G(a),\quad a>0,t\in\mathbb R.
\end{equation}

(iii) $G(a+it)$ is a continuous negative definite function of $t$ for
$a\ge 1$. 
\end{cor} 

Concerning  continuous positive and negative definite functions we
refer to
e.g. \cite{B:F}.

 Letting $a\to 0^+$
in \eqref{eq:L2}, we formally get that $1-G(it)$ is the Fourier
transform of $\varphi(s)1_{[0,\infty)}(s).$ 
This is true in the
$L^2$-sense because of Plancherel's theorem. In fact, we have

\begin{prop}\label{thm:L2} The function $\varphi$ in \eqref{eq:phi} is
  square integrable and
\begin{equation}\label{eq:Pl}
\lim_{a\to 0^+}\Int |1-G(a+it)|^2\,\frac{dt}{2\pi}=\Int
|1-G(it)|^2\,\frac{dt}{2\pi}=\int_0^\infty \varphi^2(s)\,ds.
\end{equation}
\end{prop}

The function $G$ is one-to-one when considered on the positive real
line. It is shown below that $G$ is conformal when defined in a sector
containing the positive real line. We put
$$
S(a,b)=\{ z\neq 0\, | \, a<\Arg z<b\}.
$$
\begin{prop}
\label{prop:conformal}
The function $G:S(-\pi/3,\pi/3)\to \mathbb C$ is
a conformal mapping.  
\end{prop}
 
Based on computer experiments it seems that $G$ is conformal in the
right half plane, but we have not been able to verify
this. On the other hand, $G:\mathcal A\to \mathbb C$ is not conformal.
 
\section{Proof of the properties of $G$}

In the first lemma the behaviour of $G$ close to zero and infinity is
investigated. 
\begin{lemma}\label{thm:zeroinfinity} We have
\begin{enumerate}
\item[(i)] There exist constants $A,B>0$ such that
$$
|G(z)|\le A|\Log z|+B|\Log z|^2\;\;\mbox{for}\;\; z\in\mathcal
A,|z|\le 1/2,
$$
\item[(ii)] $zG(z)\to 0$ for $z=\varepsilon e^{i\t},\varepsilon\to 0,
$  uniformly for $-\pi<\t<\pi$.
\item[(iii)] There exists a constant $C>0$ such that
$$
|1-G(z)|\le C/|z|\;\;\mbox{for}\;\; z\in\mathcal A,|z|\ge e,
$$
\item[(iv)] $G(z)\to 1$ for $z=Re^{i\t},R\to\infty,
$ uniformly for $-\pi<\t<\pi$.
\end{enumerate}
\end{lemma}

{\it Proof.} We have for $z\in\mathcal A$
$$
G(z)=\frac{z}{\Log(1+z)}\left(\Log(1+z)\Log z-(\Log z)^2\right),
$$
hence for $|z|\le 1/2$
$$
|G(z)|\le \max_{|z|\le 1/2}\left|\frac{z}{\Log(1+z)}\right|\left(|\Log
  z|\max_{|z|\le 1/2}|\Log(1+z)|+|\Log z|^2\right),
$$ 
which shows (i). 

(ii) follows from (i) since  $z(\Log z)^n\to 0$ for $n\ge 1$ and
$|z|=\varepsilon\to 0$. 

To see (iii), we note that the power series (in $1/z$)
$$
\Log(1+1/z)=\sum_{n=1}^\infty (-1)^{n-1}\frac1{nz^n},\quad |z|>1
$$
yields
\begin{equation}\label{eq:Log}
|\Log(1+1/z)|\le \sum_{n=1}^\infty 1/|z|^n =\frac{1}{|z|-1}\le
\frac{1}{e-1},\quad |z|\ge e.
\end{equation}

The power series also yields
\begin{equation}\label{eq:Log1}
z\Log(1+1/z)=1+\alpha(z)/z,\; |\alpha(z)|\le \frac{e}{2(e-1)},
\quad  |z|\ge e.
\end{equation}

Note also that $|\Log z|\ge 1$ for $z\in\mathcal A,|z|\ge e$.

Writing
$$
\frac{\Log z}{\Log(1+z)}=1+\beta(z)\Log(1+1/z),
$$
with
$$
\beta(z)=\frac{-1}{\Log(1+z)}
$$ 
we find for $z\in\mathcal A,|z|\ge e$
\begin{equation}\label{eq:Log2}
|\beta(z)|=\frac{1}{|\Log(1+z)|}\le
\frac{1}{\ln|1+z|}\le\frac{1}{\ln(|z|-1)}\le\frac{1}{\ln(e-1)}.
\end{equation}

Finally, since
$$
G(z)=(z\Log(1+1/z))\frac{\Log z}{\Log(1+z)}=(1+\alpha(z)/z)(1+\beta(z)\Log(1+1/z)),
$$
we see that
$$
z(G(z)-1)=\alpha(z)+\beta(z)z\Log(1+1/z) +\alpha(z)\beta(z)\Log(1+1/z),
$$
which by \eqref{eq:Log1} and \eqref{eq:Log2} is bounded by some constant $C>0$ for $|z|\ge e$, showing
(iii). Property (iv) follows immediately from (iii).
\hfill $\square$

\medskip
{\bf Proof of Theorem~\ref{thm:G1}}
For fixed $z\in\mathcal A$ we choose $\eps$ and $R$ such that $0<\eps<|z|<R$ and consider the positively oriented contour
$\gamma(\eps,R)$ in
$\mathcal A$ consisting of the half-circle $z=\eps e^{i\t},\t\in
[-\tfrac{\pi}2,\tfrac{\pi}2]$ and the half-lines $z=x \pm i\eps,x\le 0$
until they cut the  circle $|z|=R$, which closes the contour at the
points $-R(\varepsilon)\pm i\varepsilon$, where $0<R(\varepsilon)\to R$
for $\varepsilon\to 0$.
By Cauchy's integral theorem we have
\begin{equation}\label{eq:CIT}
G(z)=\frac{1}{2\pi i}\int_{\gamma(\eps,R)}
\frac{G(w)}{w-z}\,dw.
\end{equation}
Letting $\varepsilon$ tend to zero, the contribution corresponding to
the half-circle with radius $\varepsilon$ tends to 0 by (ii)
of Lemma~\ref{thm:zeroinfinity}.

Concerning the boundary behaviour of $G$ on the negative real line we obtain
\begin{equation*}
G(t+i0):=\lim_{\varepsilon\to 0^+}G(t+i\varepsilon)
=\left\{\begin{array}{ll}
\left(1-\frac{\ln(-t)+i\pi}{\ln|1+t|+i\pi}\right)t(\ln(-t) +i\pi),&
\;\mbox{if}\;\; t<-1\\
&\\
\left(1-\frac{\ln(-t)+i\pi}{\ln(1+t)}\right)t(\ln(-t)+i\pi),&
\;\mbox{if}\;\; -1<t<0.
\end{array}\right.
\end{equation*}
Note that $G(t+i0)$ is continuous at $t=-1$ with value $-i\pi$. Using
that $G(\overline{z})=\overline{G(z)}$, \eqref{eq:CIT} yields
\begin{equation}\label{eq:CIT1}
G(z)=\frac{1}{2\pi}\int_{-\pi}^\pi\frac{G(Re^{i\t})}{Re^{i\t}-z}Re^{i\t}\,d\t
+\frac{1}{\pi}\int_{-R}^0\frac{\Im G(t+i0)}{t-z}\,dt.
\end{equation}
In the last integral we replace $t$ by $-t$ and use that $(-1/\pi)\Im
G(-t+i0)=-\rho(t)$. Letting
$R\to\infty$ and using (iv) of Lemma~\ref{thm:zeroinfinity},
we finally get \eqref{eq:Grep}.
\hfill $\square$

\begin{rem} {\rm Feng Qi has kindly informed us about the following
    elementary proof of the observation that $1-G$ is not a Stieltjes
    function. In fact, if it were, then also $h(x)=1/(x(1-G(x))$ would
    be a Stieltjes function by the Stieltjes-Reuter-It\^o Theorem,
    cf. \cite{Re}, \cite{B1} or \cite[p.25]{B2}. In particular, $h$ will be
    decreasing, which is contradicted by the simple fact that
$1=h(1)<h(2)=1.02\cdots$.
}
\end{rem}
\medskip
{\bf Proof of Theorem~\ref{thm:G2}}
The formulas \eqref{eq:1-G} and \eqref{eq:phi} follow immediately from Theorem~\ref{thm:G1}
and it remains to prove that $\varphi$ is positive. 
Let $t_0=(\sqrt{5}-1)/2$. Then $\rho(t)>0$ for $0<t<t_0$ and
$\rho(t)<0$ for $t_0<t<\infty$ and hence
\begin{equation}\label{eq:AB}
A=\int_0^{t_0} \rho(t)\,dt>0,\quad B=\int_{t_0}^\infty\rho(t)\,dt<0.
\end{equation}
Using this notation we get 
\begin{eqnarray*}
\varphi(s)&=&\int_0^{t_0}e^{-st}\rho(t)\,dt + \int_{t_0}^\infty e^{-st}\rho(t)\,dt  \\
&\ge& \int_0^{t_0}e^{-st_0}\rho(t)\,dt + \int_{t_0}^\infty
e^{-st_0}\rho(t)\,dt=(A+B)e^{-st_0}.
\end{eqnarray*}
In the following lemma it is established that  $A+B>0$, and hence $\varphi(s)>0$ for all $s\ge 0$. \hfill $\square$

\begin{lemma}\label{thm:A>B} 
$$
\int_0^\infty \rho(t)\,dt>0.
$$
\end{lemma}
 {\it Proof.} We first establish
\begin{equation}\label{eq:int01}
\int_0^1 \rho(t)\,dt> \frac{\pi^2}{6}-\frac{1}{2}.
\end{equation}
Since  
$$
\sum_{n=0}^\infty t^n\int_0^1\binom{x}{n}\,dx=\int_0^1
(1+t)^x\,dx=\left[\frac{(1+t)^x}{\ln(1+t)}\right]_0^1=\frac{t}{\ln(1+t)}
$$
we obtain the power series expansion
\begin{equation}\label{eq:ps}
\frac{t}{\ln(1+t)}=1+\sum_{n=1}^\infty b_nt^n,\quad |t|<1;\quad
b_n=\int_0^1 \binom{x}{n}\,dx.
\end{equation}
The numbers $b_n$ are sometimes called the Cauchy numbers.
Note that for $n\ge 1$
\begin{equation}\label{eq:ps1}
0<(-1)^{n-1}b_n=\int_0^1\frac{x(1-x)\cdots(n-1-x)}{n!}\, dx\le
\frac{1}{n}\int_0^1 x\,dx=\frac{1}{2n}.
\end{equation}
By \eqref{eq:ps} we get
\begin{eqnarray*}
\int_0^1\rho(t)\,dt&=&-\frac12-2\int_0^1\ln t\,dt + 2\sum_{n=1}^\infty
(-1)^{n-1}b_n\int_0^1 t^n\ln t\,dt\\
&=&\frac32 - 2\sum_{n=1}^\infty (-1)^{n-1}\frac{b_n}{(n+1)^2}
\end{eqnarray*}
and hence using \eqref{eq:ps1}
\begin{eqnarray*}
\int_0^1\rho(t)\,dt &>& \frac32 -\sum_{n=1}^\infty \frac{1}{n(n+1)^2}
=\frac32 -\sum_{n=1}^\infty
\frac1{n+1}\left(\frac1{n}-\frac1{n+1}\right)\\
&=& \frac12 +\sum_{n=1}^\infty\frac{1}{(n+1)^2}=\frac{\pi^2}{6}-\frac12.
\end{eqnarray*}

We next show that
\begin{equation}\label{eq:int12}
\int_1^2 \rho(t)\,dt > -\frac{1}{12}-\frac{2(1+\ln 2)\ln 2}{\pi^2}, 
\end{equation}
by using the rough estimate
$$
\int_1^2 \rho(t)\,dt > -\frac{1}{\pi^2}\int_1^2 t(\ln(1-1/t))^2\,dt
$$
and
$$
\int_1^2 t(\ln(1-1/t))^2\,dt=\int_1^2 \left(t(\ln(t-1))^2+t(\ln t)^2
-2t\ln(t-1)\ln t\right)\,dt.
$$
The integral of the first two terms can be calculated because 
$$
\int t(\ln t)^2\,dt=\frac{t^2}{2}\left((\ln t)^2-\ln t+\frac12\right),
$$
and for the integral of the third term we have
\begin{eqnarray*}
2\int t\ln(t-1)\ln t\,dt &=&\left(t^2\ln t-\frac12
  t^2+\frac12\right)\ln(t-1) -\frac12 t^2\ln t +\frac12 t^2\\
&& -t\ln t+\frac32 t + \dilog(t),
\end{eqnarray*}
where
$$
\dilog(t)=\int_1^t \frac{\ln x}{1-x}\,dx.
$$ 
Since
$$
\dilog(2)=-\int_0^1 \frac{\ln(1+u)}{u}\,du=-\sum_{n=0}^\infty \int_0^1
(-1)^n\frac{u^n}{n+1}\,du=-\frac{\pi^2}{12},
$$
this leads to the expression in \eqref{eq:int12}.

We finally show
\begin{equation}\label{eq:int2infty}
\int_2^\infty \rho(t)\,dt > -\frac12.
\end{equation}
Squaring the power series for $\ln(1-u)$ yields
\begin{equation}\label{eq:ps2}
(\ln(1-u))^2=u^2\sum_{n=0}^\infty c_nu^n,\quad |u|<1;\quad
c_n=\sum_{k=0}^n\frac{1}{(k+1)(n+1-k)}.
\end{equation}
The relation $0<c_n\le 1$ for all $n$ is proved in Lemma~\ref{thm:cn}
below. Therefore, and using \eqref{eq:ps2} with $u=1/t$ it follows
that
\begin{align*}
\int_2^\infty \rho(t)\,dt &= -\int_2^\infty \sum_{n=0}^\infty
\frac{c_n}{t^{n+1}((\ln(t-1))^2+\pi^2)}\,dt\\
&> -\int_2^\infty \left(\sum_{n=0}^\infty
\frac{1}{t^{n+1}}\right)\frac{dt}{(\ln(t-1))^2+\pi^2}\\
&=-\int_2^\infty
\frac{dt}{(t-1)((\ln(t-1))^2+\pi^2)}\\
&= -\frac{1}{\pi}\int_0^\infty \frac{dx}{1+x^2}=-\frac12.
\end{align*}
Combining \eqref{eq:int01}, \eqref{eq:int12} and
\eqref{eq:int2infty} we get
$$
\int_0^\infty \rho(t)\,dt>\frac{\pi^2}{6}-\frac{1}{2}-\frac{1}{12}-\frac{2(1+\ln 2)\ln 2}{\pi^2}-\frac{1}{2}\simeq 0.3238 > 0
$$
and the lemma is proved.\hfill $\square$
\begin{rem}
A  numerical computation yields
$$
\varphi(0)=\int_0^\infty \rho(t)\,dt\simeq 0.5192.
$$
\end{rem}

\begin{lemma}\label{thm:cn} The numbers  
$$
c_n=\sum_{k=0}^n\frac{1}{(k+1)(n+1-k)},\quad n\geq 0
$$
can be written in the form 
$$
c_n=\frac{2\mathcal H_{n+1}}{n+2},
$$
where $\mathcal H_n=\sum_{k=1}^n1/k$
is the $n$'th harmonic number. Moreover,
$$
c_{n-1}-c_n=\frac{2\left(\mathcal H_n -1\right)}{(n+1)(n+2)}\ge 0,
$$
whence $1=c_0=c_1>c_2>c_3\ldots$.
\end{lemma} 
{\it Proof.} By definition we have
\begin{align*} 
c_n = \sum_{k=0}^n \int_0^1 x^k\,dx \int_0^1 y^{n-k}\,dy
&=\int_0^1\int_0^1 \frac{x^{n+1}-y^{n+1}}{x-y}\,dx\,dy\\
&= 2\int_0^1\left(\int_0^x
  \frac{x^{n+1}-y^{n+1}}{x-y}\,dy\right)\,dx\\
&=2\int_0^1
x^{n+1}\,dx \int_0^1\frac{1-t^{n+1}}{1-t}\,dt=\frac{2\mathcal H_{n+1}}{n+2}
.
\end{align*}
Using this formula we find
$$
c_{n-1}-c_n=\frac{2}{(n+1)(n+2)}\left((n+2)\mathcal H_n -
  (n+1)\mathcal H_{n+1}\right)=\frac{2\left(\mathcal H_n -1\right)}{(n+1)(n+2)}
$$
which proves the lemma. $\quad\square$

\medskip

 {\bf Proof of Corollary~\ref{thm:Bochner}}. It is well-known that if
 $F(t)$ is a continuous positive definite function on $\mathbb R$,
 then $F(0)-F(t)$ is continuous and negative definite, and a
 continuous negative definite function $H(t)$ satisfies $\Re\,H(t)\ge
 H(0)\ge 0$, see \cite{B:F}. Therefore (ii) follows from (i), and (iii)
 follows from (ii) because $G(a)\ge 0$ for $a\ge 1$. $\quad\square$

\medskip

 {\bf Proof of Proposition~\ref{thm:L2}}. By (i) and (iii) of Lemma~\ref{thm:zeroinfinity} it is
 clear that 
$$
\Int|1-G(it)|^2\,\frac{dt}{2\pi}<\infty
$$
and that dominated convergence can be applied to yield the first
equality in \eqref{eq:Pl}. By Plancherel's theorem $1-G(it)$ must be
the Fourier transform of a square integrable function, which is the
$L^2$-limit of \eqref{eq:plan}, hence equal to
$\varphi(s)1_{[0,\infty)}(s).$
\hfill $\quad\square$  

Before proving Proposition \ref{prop:conformal} we give Lemma
\ref{lemma:G-prime}.
\begin{lemma}
\label{lemma:G-prime}
For $z\in S(0,\pi/3)$ we have $
\Im G'(z)<0$.
\end{lemma}
{\it Proof.} From the relation \eqref{eq:Grep} it follows that
\begin{equation}
\label{eq:G-prime-rep}
\Im G'(re^{i\theta})=
-2r\sin\theta\int_0^{\infty}\frac{r\cos \theta+t}{((r\cos
  \theta+t)^2+(r\sin \theta)^2)^2}\rho(t)\, dt.
\end{equation}
We claim that for fixed $r>0$ and $\theta \in [0,\pi/3]$ the function 
$$
k(t)=\frac{r\cos \theta+t}{((r\cos
  \theta+t)^2+(r\sin \theta)^2)^2}
$$
is decreasing. Indeed, it follows that 
$$
k'(t)=\frac{(r\sin \theta)^2-3(r\cos \theta+t)^2}{((r\cos
  \theta+t)^2+(r\sin \theta)^2)^3}
$$
and the numerator is negative for all $t>0$ because $\sin^2
\theta\leq 3\cos ^2 \theta$ for $\theta\in [0,\pi/3]$. 

This implies
\begin{align*}
\int_0^{\infty}\frac{r\cos \theta+t}{((r\cos
  \theta+t)^2+(r\sin \theta)^2)^2}\rho(t)\, dt
&=
\int_0^{t_0}k(t)\rho(t)\, dt+
\int_{t_0}^{\infty}k(t)\rho(t)\, dt\\
&\geq k(t_0)\left(\int_0^{t_0}\rho(t)\, dt+
\int_{t_0}^{\infty}\rho(t)\, dt\right),
\end{align*}
where $t_0=(\sqrt{5}-1)/2$. From Lemma \ref{thm:A>B} it follows that the integral above is positive. From \eqref{eq:G-prime-rep} we now obtain that
$$
\Im G'(re^{i\theta})=
-2r\sin\theta\int_0^{\infty}\frac{r\cos \theta+t}{((r\cos
  \theta+t)^2+(r\sin \theta)^2)^2}\rho(t)\, dt<0.
$$
This proves the lemma.\hfill $\square$

{\bf Proof of Proposition \ref{prop:conformal}.} From \eqref{eq:Grep}
it follows that  
$$
\Im G(x+iy)=y\int_0^{\infty}\frac{\rho(t)}{(x+t)^2+y^2}\, dt.
$$
Here $t\mapsto 1/((x+t)^2+y^2)$ is a decreasing function of $t$ and it
follows as in Lemma \ref{lemma:G-prime} that $\Im G(x+iy)>0$ for
$x>0$ and $y>0$ and also  $\Im G(x+iy)<0$ for
$x>0$ and $y<0$. Hence it is enough to show that $G$ is one-to-one in
the sector $S(0,\pi/3)$.
 
 For $z_1$ and
$z_2$ belonging to the sector $S(0,\pi/3)$ we have
$$
G(z_2)-G(z_1)= \int_{\gamma(z_1,z_2)}G'(w)\, dw,
$$
where $\gamma(z_1,z_2)$ is the straight line segment from $z_1$ to $z_2$. Thus
$$
G(z_2)-G(z_1)= (z_2-z_1)\int_0^1G'(z_1+t(z_2-z_1))\, dt\neq 0,
$$
when $z_1\neq z_2$ since $\Im G'(w)<0$ for $w\in S(0,\pi/3)$ by Lemma \ref{lemma:G-prime}. This shows that $G$ is one-to-one in $S(0,\pi/3)$.\hfill $\square$

\end{document}